\documentstyle{article}
\renewcommand{\baselinestretch}{1.3}

\newtheorem {th}{Theorem}
\newtheorem {lem}[th]{Lemma}

\def\Cox{\hfill \Box}

\def\deq{\, {\stackrel {def} {=}}}

\def\sf{\sigma\mbox{-field}}
\def\dd{\delta}

\def\ee{\epsilon}
\def\E{{\bf{E}}}
\def\P{{\bf{P}}}
\def\N{\hbox{I\kern-.2em\hbox{N}}}

\def\nbyn{[n \times n]}

\def\U2{(U \times U)}
\def\S{{\cal{S}}}
\def\R{{\cal{R}}}

\def\F{{\cal{F}}}
\def\cE{{\cal{E}}}
\def\sgn{\mbox{sign}}
\def\|{\, | \, }
\def\one{{\bf 1}}
\def\0{\hat{0}}
\def\1{\hat{1}}

\begin{document}

\begin{titlepage}

\begin{center}
{\large \bf A SHUFFLE THAT MIXES SETS OF ANY FIXED SIZE MUCH
  FASTER THAN IT MIXES THE WHOLE DECK} \\
\end{center}
\vspace{5ex}
\begin{flushright}
Robin Pemantle \footnote{Research supported in part by a
National Science Foundation Grant \# DMS 9103738} \\

\end{flushright}

\vfill

{\bf ABSTRACT:} \break
Consider an $n$ by $n$ array of cards shuffled in the following
manner.  An element $x$ of the array is chosen uniformly at random;
Then with probability $1/2$ the rectangle of cards above and to
the left of $x$ is rotated 180 degrees, and with probability $1/2$
the rectangle of cards below and to the right of $x$ is rotated
180 degrees.  It is shown by an eigenvalue method
that the time required to approach the uniform
distribution is between $n^2/2$ and $c n^2 \ln n$ for some constant
$c$.  On the other hand, for any $k$ it is shown that the
time needed to uniformly distribute a set of
cards of size $k$ is at most $c(k) n$, where $c(k)$ is
a constant times $k^3 \ln (k)^2$.  This is established via
coupling; no attempt is made to get a good constant. 
\vfill

\noindent{Keywords:} shuffle, array, randomization time, 
coupling, eigenvalue

\noindent{Subject classification: } 60B15 , 60J15

\end{titlepage}

\section{Introduction}

Consider $n^2$ playing cards, numbered $1 , \ldots , n^2$, in an
$n \times n$ array; the set of positions in this array is denoted by
$$\nbyn = \{ (i,j) : 1 \leq i,j \leq n \} ,$$
with $(1,1)$ in the upper-left corner.
For $1 \leq i,j \leq n$, let $\pi_{ij}$ be the
permutation that sends the card in the $(r,s)$ position to
the $(i+1-r,j+1-s)$ position if $r \leq i$ and $s \leq j$, and
does not change the position of the card otherwise.  In other
words the rectangle of size $i \times j$ in the upper-left corner
gets rotated by $180^\circ$ and the remaining cards are unmoved.
(The $(1,1)$ position is in the upper left, following matrix rather
than Cartesian notation.)
Let $\pi_{ij}'$ denote the shuffle that does the same for the
lower right corner, so that the card in the $(r,s)$ position 
is moved to position $(n+i-r,n+j-s)$ if $r \geq i$ and $s \geq j$
and is otherwise unmoved.  Questions about how rapidly this type of
permutation mixes an array were inspired by a Macintosh screensaver.

Suppose first that the cards are shuffled by waiting a mean one exponential
amount of time, then picking $i$ and $j$ uniformly at random and
performing the shuffle $\pi_{ij}$.  (Setting the problem in continuous 
time avoids the later use of more complicated versions of theorems 
in~\cite{Di} and~\cite{DS} that take parity into account.)
After time $t$, the resulting distribution 
$\S_0^t$ on permutations of the $n^2$ positions is given by
$$\S_0^t = \exp (t  (\S_0 - 1)) \deq \sum e^{-t} {t^k \over k!} 
   \S_0^{(k)}$$
where $\S_0^{(k)}$ is the $k$-fold convolution of the measure
$\S_0 = n^{-2} \sum_{i,j = 1}^n \dd_{\pi_{ij}}$.  Here and throughout,
random walks on the space of card configurations are identified 
with random walks on the symmetric group; in particular, when
discussing two coupled shuffles, it will be convenient to be able
to refer to the positions $\sigma (x)$ and $\tau (x)$ of the same 
card $x$ in two arrays starting from two arbitrary configurations, 
one permuted by $\sigma$ and the other by $\tau$.  

The card in the $(n,n)$ position is unlikely to move before time 
$cn^2$, which gives an easy lower bound on the time needed to
randomize the layout.  More precisely, if $A$ is the set of
permutations fixing $(n,n)$, then $\S_0^t (A) \geq e^{-t/n^2}$
since $e^{-t/n^2}$ is the probability that the card in position 
$(n,n)$ is never moved at all.  Thus 
$$|S_0^t - U| \geq e^{-t/n^2} - {1 \over n^2} , $$
where $U$ is the uniform measure and $| \cdot |$ is the total
variation distance.  When $t << n^2$, therefore, the total variation
distance is near one and the deck is not well shuffled.
The same lower bound may be obtained by counting:
the total number of permutations of $n^2$ cards is
$$n^2 ! = \exp ((2+o(1)) n^2 \log n) ,$$
whereas the set $A_k$ of permutations reachable in $k$ shuffles
is at most $n^{2k}$.   Thus, letting $k = \lfloor (1 + \ee) t \rfloor$, 
\begin{eqnarray*}
|S_0^t - U| & \geq & |\S_0^t (A_k) - U(A_k)| \\[2ex]
& = & 1 + o(1) - \exp [ 2 \log n (k - (1+o(1))n^2)] ,
\end{eqnarray*}
which is near 1 when $t << n^2$.  It will be seen (Theorem~\ref{tv time} 
below) that the time to randomization is at most a constant 
times $n^2 \ln (n)$.

The shuffle becomes more interesting if permutations $\pi'_{ij}$
are also allowed.  If each $\pi_{ij}$ and $\pi'_{ij}$ occurs at rate
$1/(2n^2)$, the distribution resulting at time $t$ will be
$$\S^t \deq \exp (t (\S - 1)) \deq \sum e^{-t} {t^k \over k!} \S^{(k)}$$
where $\S$ gives probability $1/(2n^2)$ to each $\pi_{ij}$ and
to each $\pi'_{ij}$.
(The dependence of $\S$ and $\S_0$ on $n$ is suppressed in the notation.)
Now the cards that take the longest to move are in positions $(1,n)$ 
and $(n,1)$ and these will each be moved by time $cn$ with probability
$1 - e^{-c/2}$.  Thus the first argument above shows only that the
deck is not at all shuffled by time $t << n$.  The counting argument
from before does better: setting $k = \lfloor (1+\ee) t \rfloor$ shows that
$|\S^t - U| \approx 1$ when $t << n^2$.  On the other hand, 
it will be shown that
the positions of any set of cards of any fixed size, $k$, will 
be jointly randomized by time $cn$ as $n \rightarrow \infty$.  
(By altering the shuffle again so that it may choose rectangles
in the lower left and upper right corners as well, this time can
be reduced to a constant when $k=1$, but not for $k \geq 2$, since
a pair of neighboring cards will always be stuck together for expected 
time $cn$.)
This is the only shuffle I know of with the property that the
time to randomization differs from the time to randomize subsets
of any fixed size by factors greater than $\mbox{poly-log }(n)$.
In fact, $k$ may be allowed to increase with $n$, in such a way that
the time to randomize any $k$ cards is still much less than 
the time to total randomization.  To quantify this, say that an event
$A$ is measurable with respect to cards $x_1 , \ldots , x_k$
if $A$ is a set of permutations of the form
$\{ \pi : (\pi (x_1) , \ldots , \pi (x_k)) \in B \}$ for some
cards $x_1 , \ldots , x_k$, where $B$ is a subset of $k$-tuples
of distinct positions in the array $\nbyn$.
Define the $k$-set distance to uniformity of a distribution $\R$, 
denoted $||\R - U||_k$, to be $\sup_A \R(A) - U(A)$ as $A$ ranges
over events measurable with respect to the positions of some
set of $k$ cards; setting $k=n$ recovers the total variation distance.  

\begin{th} \label{k-sets}
There exists a constant $c$ such that for any $n$ and any $k$ with
$1 < k < n$,
$|| \S^t - U ||_k < 1/2^j$ whenever $t > c k^3 (\ln(k))^2 nj$.  
\end{th}

\begin{th} \label{tv time}
For any $\ee  > 0$, $\lim_n | \S^t - U | =1$ when $t = (1-\ee) n^2/2$.
On the other hand there is a constant $c$ for which $| \S^t - U| < 1/2^j$
whenever $t > cj n^2 \ln (n)$.  The same is true with $\S$ replaced by
$\S_0$.
\end{th}

The author wishes to thank Martin Hildebrand for helpful comments
toward the revised draft of this manuscript.
The proofs of both theorems are based on techniques developed by
Diaconis and others~\cite{Di,DS}.  In particular, the second part
of Theorem~\ref{tv time} uses eigenvalue machinery (the first
part is just a counting argument) and the proof of Theorem~\ref{k-sets}
is a coupling argument.  No new theory is developed in this paper, rather
it is hoped that the example is interesting.  

\section{Proof of Theorem \protect{\ref{k-sets}}}

Theorem\ref{k-sets} is proved via a series of lemmas that establish it
for small values of $k$. Do not count on an unsubscripted $c$ to denote 
the same quantity from line to line.

\begin{lem} \label{k=1}
There exists a constant $c$ such that for any $n$,
$|| \S^t - U ||_1 < 1/2^j$ whenever $t > cjn$.  
The author wishes to thank Martin Hildebrand for helpful comments
toward the revised draft of this manuscript.
\end{lem}

\begin{lem} \label{k=2}
There exists a constant $c$ such that for any $n$,
$|| \S^t - U ||_2 < 1/2^j$ whenever $t > cj n$.  
\end{lem}

\begin{lem} \label{k=3}
There exists a constant $c$ such that for any $n$,
$|| \S^t - U ||_3 < 1/2^j$ whenever $t > cj n$.  
\end{lem}

To get from each lemma to the next, and thence to the theorem, the
following type of coupling argument is used.  For each finite set of cards
$(x_1 , \ldots , x_k)$, a Markov chain $\{ (\sigma_t , \tau_t ) : t \geq 0 \}$ 
is defined on pairs of permutations of $n^2$ cards.  It is a coupling of two
copies of the shuffle $\S$ in the sense that the marginal on either
coordinate is Markov with transitions from $\sigma$ to $\sigma \pi_{ij}$
or $\sigma \pi'_{ij}$ at rates $1/(2n^2)$ each, and that from some point
onward $\sigma (x_i)$ will equal $\tau (x_i)$ for all $i$.  (At this time the 
coupling is said to have succeeded, the initial configurations of cards having
been any two arbitrary configurations.)  Furthermore, there are
constants $c', \dd > 0$ independent of the cards $x_1 , \ldots , x_k$
such that for any pair $(\sigma_0 , \tau_0)$, the probability that 
the coupling will succeed by time $c' n$ is at least $\dd$.  Repeating this
coupling $j \lceil \log (1/2) / \log (1-\dd) \rceil$ times and letting 
$c = c' \lceil \log (1/2) / \log (1-\dd) \rceil$ gives a coupling for which 
the probability that $\sigma_t (x_i) = \tau_t (x_i)$ for all $t \geq cjn$ 
and $1 \leq i \leq k$ is at least $1-1/2^j$.  Since $x_1 , \ldots , x_k$
were arbitrary as were the two initial configurations, this implies the desired
conclusion.  It remains to exhibit the couplings, which will be done in the
notation of this paragraph and without any thrift in choices of constants.
To avoid drowning in a mire of greatest-integer brackets, ignore them, 
i.e., assume without loss of generality that $n$ is divisible by all of
the integer constants that arise in the proofs.  Also, names such as $A$
and $B$ will be assigned anew for each lemma.

\noindent{Proof} of Lemma~\ref{k=1}:  
For each starting position $(i,j) \in \nbyn$,
consider the set of possible positions to which 
a card in that position may jump under a single permutation, 
$\pi_{rs}$ or $\pi_{rs}'$.  This
is just the set $\{ (a,b) : (n+1-a-i)(n+1-b-j) \geq 0 \}$; 
pictorially, rotate by $180^\circ$ to get the point $(n+1-i,n+1-j)$, 
then divide the array into (unequal) quadrants meeting there
and the possible jump set will consist of the upper-left and lower-right
quadrants; the jump set is the shaded region in figure 1.  
Let $A \subseteq \nbyn$
be the region $i,j \leq n/3$ and let $B$ be the region $i,j \geq 2n/3$;
see figure 2.  Observe that for any card $x_1$, in any position $(i,j)$,
the rate at which $x_1$ jumps into the region $A \cup B$ is at least 
$1/(3n)$.  Indeed, the area of intersection of $A \cup B$ with the
shaded region in figure 1 is minimized when $(i,j) = (1,n)$ or
$(i,j) = (n,1)$.  It is therefore possible
to construct a coupling where at rate $1/3n$, independent of the past,
both coordinates, $\sigma$ and $\tau$, 
simultaneously jump to permutations for which the card $x_1$ is
in $A \cup B$.  Call the first time this happens $T$.
From the pictorial description of the jump set, it follows 
that any two positions in $A \cup B$
have at least $n^2/3$ positions in common to which both may jump
($n^2 / 9$ suffices for our argument and is more immediate).

To finish the argument, let $C$ denote the set of positions 
reachable in a single jump from both $\sigma_T (x_1)$ and 
$\tau_T (x_1)$.  Then the probability that the process
$\{ \sigma_t (x_1) : T < t \leq T+1 \}$ contains precisely
one jump and that $\sigma_{T+1} (x_1) \in C$ is at least
$|C|/2n^2$ times the probability of exactly one jump,
and therefore at least $(1/6) e^{-1}$.  The same is true
for the process $\{ \tau_T (x_1) : T < t \leq T+1 \}$.
Thus the laws of $\sigma_{T+1} (x_1)$ and $\tau_{T+1} (x_1)$
both dominate a measure uniform on $C$ with total mass $e^{-1}/6$,
and the coup[ling may be extended to time $T+1$ in such a way
that the $\P (\sigma_{T+1} (x_1) = \tau_{T+1} (x_1)) \geq 
e^{-1} / 6$.  The coupling then succeeds
in time $3n +1$ with probability at least $\P (T \leq 3n) e^{-1} /6
\geq e^{-1} (1-e^{-1}) /6$ which proves the lemma.      $\Cox$

\noindent{Proof} of Lemma~\ref{k=2}:  A useful observation is that if
cards $x_1$ and $x_2$ are both some minimal distance $d$ from any
edge of the array, and some permutation $\pi_{ij}$ is applied which 
moves $x_1$ but not $x_2$, then further application of any $\pi_{kl}$ 
with $n-d/2 \leq k,l \leq n-d/4$ sends both cards to positions at least
$d/4$ distant from any edge of the array.
Some notation for distance from the set of positions distant from any
edge will also be useful.  Let $A_j \subseteq \nbyn$
be the set of positions 
$$\{ (i,k) : {n \over 5 \cdot 2^j} \leq i,k \leq n - 
   {n \over 5 \cdot 2^j} \}.$$ 
Let $B \subset \nbyn^2$ denote the set
$$\{ ((i_1 , j_1) , (i_2 , j_2)) \in
   (A_4)^2 : \max (|i_1 - i_2| , |j_1 - j_2|) \geq n/40 \ ;$$ 
of pairs of positions in $A_4$ separated by at least $n/40$ in at 
least one coordinate.  Define $B_0$ to be the set of pairs of positions, 
one of which is in $A_2$ and the other of which has both coordinates
less than $n/40$.  Let $C \subset \nbyn^2$ be the set 
$$\{ ((i_1 , j_1) , (i_2 , j_2)) \subset (A_6)^2 : 
   \min (|i_1 - i_2| , |j_1 - j_2|) \geq n/160 \} .$$ 
Finally, let $D$ be the set of pairs of coordinates
$\{ ((i_1 , j_1) , (i_2,j_2)) : i_1 , j_1 < n/3, i_2 , j_2 > 2n/3 \}$.

Pick any distinct cards $x_1$ and $x_2$, and suppose the positions,
$(i_1, j_1)$ and $(i_2 , j_2)$ of both cards are in $A_2$.  
Either $i_1 \neq i_2$ or $j_1 \neq j_2$; assume 
without loss of generality that $i_1 \neq i_2$, since the
argument is symmetric in $i$ and $j$; furthermore, assume
without loss of generality that $i_1 < i_2$, since the argument
is symmetric in the two copies of the shuffle.  
If we choose $j$ so that
$j_1 \leq j \leq j_1 + n/20$, then the permutation
$\pi_{i_1 j}$ moves $x_1$ to a position $(1,b)$ with $b \leq n/40$
and does not move $x_2$.  The positions of the cards now differ by 
at least $n/40$ in the second coordinate.  Since any permutation
$\pi_{kl}$ with $k,l > 39n/40$ will move both cards, it will
also preserve their separation; applying the 
observation at the beginning of this proof (with
$d = n/20$) shows that
there are at least $n^2 / 6400$ permutations $\pi_{ij}$ whose 
further application will result in the cards $x_1$ and $x_2$ having
a pair of positions in $B$.  It has thus been shown that
\begin{quote}
Whenever $x_1, x_2 \in A_2$, the rate of jumping to a pair of
positions in $B_0$ is at least $1/(160n)$; when the pair of
positions is in $B_0$ then the rate of jumping to a pair in $B$
is at least $1/12800$.  
\end{quote}
Similar reasoning shows that whenever the pair of positions of 
$x_1$ and $x_2$ is in $B$, the probability that the pair of positions
will be in $C$ two jumps later is at least a constant, $c$:
there are at least $n^2/25600$ permutations $\pi_{ab}$ moving
one card into the the region 
$$\{ (r,s) : 1 \leq r,s \leq n/160 \}$$
while keeping the other card fixed; these also separate the
cards by at least $n/80$ in both coordinates; from here,
any $\pi_{rs}$ with $319n/320 \geq r,s \geq 159n/160$ will land
the pair of positions of $x_1$ and $x_2$ in $C$.  

A final observation along these lines is that whenever the pair of 
positions of $x_1$ and $x_2$ is in $C$, 
the probability of finding the pair in $D$ three jumps 
later is at least another constant.  The three moves which 
may be necessary are: if $x_2$ is above and to the left of
$x_1$, then apply any $\pi_{ij}$ with $i,j \geq 159n/160$ 
(otherwise, skip this step); 
now if $(i_1, j_i)$ is the new position of $x_1$, then
$i_1, j_1 \leq 319n/360$ and any $\pi_{kl}$ with $(i_1,j_1) 
\leq (k,l) \leq (i_1, j_1) + (n/320,n/320)$ will move $x_1$ into the 
upper-left corner without disturbing $x_2$; the separation 
between the cards is still at least $n/320$ in at least one 
coordinate, and the coordinates $(i_2,j_2)$ of the second
card are at least $n/320$, so there are at least $n^2/102400$
$\pi_{kl}'$ moves that will get $x_2$ into the lower-right
corner without disturbing $x_1$.

A useful and self-evident principle when coupling two identical 
copies of a countable recurrent Markov chain
is that if the rate to jump from each state in the set $\Theta$
into the set $\Xi$ is at least $\dd$, then a coupling $\{X_t , Y_t \}$
and a time $T$ exist such that $X_{T-}, Y_{T-} \in \Theta$,
$X_T , Y_T \in \Xi$, and such that the Lebesgue measure of
$\{ t < T : X_t , Y_t \in \Theta \}$ has exponential distribution
with mean $1/\dd$.  [One way to establish this is to define two
independent copies $\{X_t' , Y_t'\}$, altered in any way that reduces
the jump rate into $\Xi$ by $\dd$ at each state in $\Theta$, to
let $Z$ be an independent poisson process of rate $\dd$, to
let $T$ be the first time $t$ at which $Z_{t-} \neq Z_t$ 
while $X_t , Y_t \in \Theta$, and to let $X_t = X_t'$ and
$Y_t = Y_t'$ for $t' < t$, while $X_T$ and $Y_T$ jump into
$\Xi$ with whatever distribution was subtracted before, and
then the two evolve independently.] 

Thus the lower bound on 
the rate of jumping from a pair in $A_2$ to the set $B_0$ 
gives rise via this principle to a coupling $\{\sigma_t , \tau_t \}$
and a time $T$ at which $\sigma$ and $\tau$ simultaneously 
jump into $B_0$.  Use this coupling just up to the time $T$, and
then for $T < t < T+6$, let $\sigma$ and $\tau$ evolve 
independently.  Now essentially copy the argument at the end
of the proof of Lemma~\ref{k=1}.  The probability of precisely 6 jumps
occurring in $\sigma_t$ in the interval $(T,T+6]$ is $e^{-6} 6^6/6!
>1/7$; conditional on this, the probability that the pair
of positions of $x_1$ an $x_2$ under $\sigma_{T+6}$ is in $D$
is at least the product of the three constants above (one constant to get
to $B$ in one jump, one to get to $C$ in two more jumps and one to get
to $D$ three jumps after that).  Since $\tau_t$ behaves identically, the
probability of the event $G$ is at least a constant, where $G$ is the
event that the pairs of positions of $x_1$ and $x_2$ 
under both $\sigma_{T+6}$ and $\tau_{T+6}$ are in $D$.

Finally, observe that conditional on $G$, $\sigma$ and $\tau$
may be coupled by time $T+8$ with probability bounded away from zero:
let $\sigma$ and $\tau$ both jump exactly twice, using some 
$\pi_{i_1 ,j_1}$ and $\pi_{i_2,j_2}$ (as in the proof of the
preceding lemma) to send $x_1$ to the 
same position in $[n/6 \times n/6]$ and using some 
$\pi_{i_1 ,j_1}'$ and $\pi_{i_2,j_2}'$ to send $x_2$ to
the same position in the lower-right square of this size.  All that
remains is to bound the stopping time, $T$.

By the previous lemma there is a $k$ such that $t > kn$ implies 
$||\S^t - U||_1 < .01$.  This implies that for $t > kn$ and any card 
$x$, $\P (\S^t (x) \in A_2) \geq U (A_2) - .01 = .8$.  Thus the two 
independent copies of the Markov chain $\{ \sigma_t' \}$ and
$\{ \tau_t' \}$ used to construct the coupling must satisfy
$$\P (\sigma_t' (x_1) , \sigma_t' (x_2) , \tau_t' (x_1) , \tau_t' (x_2)
  \in A_2) \geq 1 - 4(1-.8) = .2 $$
for any $t > kn$.  In particular this implies that if $M \subseteq [kn , 2kn]$
is the set of times $t$ for which the positions of $\sigma_t' (x_1) , 
\sigma_t' (x_2) , \tau_t' (x_1)$ and $\tau_t' (x_2)$ are all in $A$, then 
$$.2 kn \leq \E \lambda (M) \leq .1 kn + n \P (\lambda (M) > .1 kn) , $$
where $\lambda$ is Lebesgue measure, and solving this gives
$\P (\lambda > .1 k n) \geq .1$.  The coupling is constructed so that 
$$\P (T < 2kn \| \lambda (M)) \geq 1 - \exp (-\lambda (M) / 160n) .$$
Thus $\P (T < 2kn) \geq (.1) (1 - \exp (-k/1600))$.  

This, together with the success of the coupling by time $T+8$ 
with constant probability, proves that the coupling succeeds
by time $2kn + 8$ with some constant probability, which
suffices to prove the lemma, since the coupling may be restarted
at times that are multiplies of $2kn+8$ until is succeeds.  $\Cox$

\noindent{Proof} of Lemma~\ref{k=3}:  This proof uses similar moves
to the last proof, so only the new part will be described.  Let $x_1, x_2$
and $x_3$ be any three cards.  By the previous lemma, choose a
$k$ for which $||\S^t - U||_2 \leq 1/4$ when $t \geq kn$.  Construct the
coupling by first letting $\sigma$ and $\tau$ evolve independently for
time $kn$.  Let $(a_j^1 , a_j^2)$ denote the position of $\sigma_t (x_j)$
and $(b_j^1 , b_j^2)$ denote the position of $\tau_t (x_j)$; for convenience,
define $a_0^1 = b_0^1 = a_0^2 = b_0^2 = 1$ 
and $a_4^1 = b_4^1 = a_4^1 = b_4^2 = n$.  Let
$$M_t = \min \{ |a_i^k - a_j^k| ,  |b_i^k - b_j^k| , |a_i^k - b_j^k|,
 : k = 1,2; i \neq j ; 0 \leq i,j \leq 4 \}.$$
Thus under both $\sigma$ and $\tau$, all cards $x_1 , x_2$ and $x_3$ 
are separated in each coordinate by $M_t$ from each other and from 
the boundary of the array, and for $i \neq j$, $\sigma_t (x_i)$ and
$\tau_t (x_j)$ are separated as well.  

Under the product uniform distribution, $\U2$, observe
$${\U2} (M_t \leq n/240) \leq .3 ;$$
this is because the event $\{ M_t \leq n/240 \}$ is the union of
36 events of probability at most $1/120$: 12 events that some
coordinate of some card under one of $\sigma_t$ or $\tau_t$ is
within $n/240$ of 0 or $n$, 12 events that some coordinate of
$\sigma_t (x_i)$ is too close to the same coordinate of 
$\tau_t (x_j)$, 6 events that some $\sigma_t (x_i)$ and 
$\sigma_t (x_j)$ are within $n/240$ in some coordinate, and  
6 events that some $\tau_t (x_i)$ and $\tau_t (x_j)$ 
are within $n/240$ in some coordinate.

Therefore $\P (M_{kn} \leq n/240) \leq 3/4$, by choice of $k$, 
since $M_{kn}$ is an event depending only on the positions of two cards.
Conditional on $M_{kn} > n/240$, $\sigma_{kn+5}$ and
$\tau_{kn+5}$ may be coupled so that the positions of all three cards 
$x_1, x_2$ and $x_3$ are the same under $\sigma$ and $\tau$ with 
probability bounded away from zero.  The five moves that may be
necessary are: (1) couple $\sigma (x_1)$ and $\tau (x_1)$ by moving
them both to the upper left $n/720 \times n/720$ square; (2) move this 
coupled card into the bottom right $n/1440 \times n/1440$ square
by time $kn+2$; (3) couple $\sigma (x_2)$
and $\tau (x_2)$ in an even smaller upper-left region; (4) move
$x_2$ to the region in the lower-right (but not all the way in the corner)
defined by $\{ (i,j) : n/360 < i,j < n/720 \}$; note that this does not
disturb $x_1$; (5) couple $x_3$.     $\Cox$

\noindent{Proof of} Theorem~\ref{k-sets} from Lemma~\ref{k=3}:
The method used to prove Lemma~\ref{k=3} may be generalized
to any $k$ but the coupling time is then exponential in $k$.  To get
a power law in $k$, it is necessary to construct a less wasteful coupling.
When $k \geq \sqrt{n}$, $k^3 n > n^2 \ln n$, and Theorem~\ref{k-sets}
is subsumed in Theorem~\ref{tv time}.  So no generality is lost in assuming
that $k < \sqrt{n}$.  Fix any $k$ cards, $x_1 , \ldots , x_k$.
A sequence of stopping times will be defined at which the probabilities
of certain ``good'' events occurring in the near future is large.  The
stopping times are called $\{ T(u,v) : 1 \leq u \leq k , 1 \leq v \leq l(u) \}$
and $\{ T_j : 0 \leq j \leq k \}$ and when $j \geq 1$, they satisfy
$$T_{j-1} < T(j,1) < T(j,1)+1 \leq T(j,2) < \cdots \leq T(j,l(j))
  < T(j,l(j)) + 1 = T_j .$$
Informally, at each $T (u,v)$, either something good happens one time unit
later, in which case $T_u = T(u,v) + 1$ and $l(u) = v$, or else we wait
for the next auspicious time, $T(u,v+1)$.  

Describing the behavior of the coupling between times $T(u,v)$ and 
$T(u,v)+1$ takes a little notation, but at all other times the construction
is simple.  Let $(\sigma_t , \tau_t)$ evolve independently until time
$T_0$.  For $t \in [T_j , T (j+1 , 1)]$ and for $t \in [T(u,v)+1, T(u,v+1)],
v < l(u)$, let $\sigma$ and $\tau$ evolve in parallel, so that 
$\sigma$ jumps to $\sigma \pi$ if and only if $\tau$ jumps to $\tau \pi$.
No technical problems arise in switching between these behaviors as
long as the $T(u,v)$ are honest stopping times and the event $\{l(u) = v \}$
is in the $\sf$ of events up to time $T(u,v)+1$.  

To handle the remaining times, define $W(t)$ to be the set 
$\{ s \leq k : \sigma_t(x_s) = \tau_t (x_s) \}$.  Informally,
this is the set of cards whose positions are the same under $\sigma$
and $\tau$ at time $t$.  Since the coupling depends on knowing
something about the configurations at times $T(u,v)$, we begin by 
defining those.  First, define
$$T_0 = \inf \{ t \geq 0 : \sigma_t (x_s) \neq \tau_t (x_{s'})
    \mbox{ for all } s ,s' \leq k \} .$$
Clearly this is a stopping time, and $W(T_0) = \emptyset$.  It will be
verified inductively that 
\begin{equation} \label{W size}
W (s) \subseteq W(t) \mbox{ for } T_0 \leq s \leq t, \; |W(T_j)| = j, \mbox{  and }
    |W(T(u,v))| = u-1.
\end{equation}
It will also be verified that $\sigma_t (x_s) \neq \tau_t (x_{s'})$ for all
$t \geq T_0$ and $s \neq s'$.  Since $\sigma$ and $\tau$ move in parallel
except on $t \in [T(u,v) , T(u,v)+1]$ and since these two statements are
true at time $T_0$, we need only verify that they remain true over the
time intervals $[T(u,v) , T(u,v) + 1]$.
For any $u \leq k$ and $1 < v \leq l(u)$, define
$$T (u,v) = \inf \{ t \geq T(u,v-1) + 1 : \begin{array}{l}
  \exists s = s(u,v) \notin W(T_{u-1}) \mbox{s.t.} \sigma_t (x_s) , 
   \tau_t (x_s) \in [1,{ n \over 3\sqrt{k}}]  \times [1,{ n \over 3\sqrt{k}}] \\
\mbox{ and } \sigma_t (x_{s'}) , \tau_t (x_{s'}) \notin [1,{ n \over 2\sqrt{k}}]
    \times [1,{ n \over 2\sqrt{k}}] \mbox{ for }  s' \neq s \} . \end{array}$$
Define $T(u,1)$ identically, but with $T_{u-1}$ in place of $T(u,v-1)+1$.
Informally, $T(u,v)$ is the first time after $T(u,v-1)+1$ (or $T_{u-1}$
if $v=1$) that some card $x_{s'}$ not yet in $W$ is sent to a square region 
in the the top-left corner by both $\sigma$ and $\tau$, while all other cards
are sent to a region in the lower-right that is the complement of a slightly
larger square region.  Clearly, these are stopping times and $W$ cannot
change on $[T(u,v-1) + 1 , T(u,v)]$ because $\sigma$ and $\tau$ are 
evolving in parallel.  

For each $a,b \leq n/(6\sqrt{k})$, there are unique $i(a,b), j(a,b)
\leq n / (2\sqrt{k})$ for
which $\pi _{ij} [\sigma_{T(u,v)} (x_s)] = (a,b)$, 
while $\pi_{ij} [\sigma_{T(u,v)} (x_{s'})] = \sigma_{T(u,v)} (x_{s'} )$ 
for $s' \neq s$.  The same is true with $\sigma_{T(u,v)}$ replaced by 
$\tau_{T(u,v)}$; call these $i^* (a,b)$ and $j^*(a,b)$.  It is therefore
possible to choose a pair $(\pi , \pi^*)$ in such a way that
each of $\pi$ and $\pi^*$ is uniform over $\{ \pi_{xy} : 1 \leq x,y 
\leq n \}$, that 
\begin{equation} \label{yes prob}
\P (\pi = \pi_{i(a,b) , j(a,b)} , \pi* = \pi_{i^* (a,b),
    j^* (a,b)}) \geq {1 \over 36 k} \; , \end{equation}
and that with probability one, either $\pi = \pi^*$ or else 
\begin{equation} \label{stayeq}
\pi = \pi_{xy} , \pi^* = \pi_{x^* y^*} \mbox{ for some } x , y , x^* , y^*
   \leq {1 \over 2 \sqrt{k} }  .
\end{equation}
For a single shuffle, $\S^t$, the probability of precisely one jump occurring
in a unit of time and that jump being a $\pi_{ij}$ rather than a $\pi_{ij}'$
is $1/(2e)$.  By this observation and~(\ref{yes prob}) and~(\ref{stayeq}),
we may construct the coupling for $t \in [T(u,v) , T(u,v)+1)]$ so that
with probability $1 - 1/(2e)$ the two processes $\sigma$ and $\tau$
evolve in parallel, jumping either zero times, more than once, or jumping
exactly once by some $\pi_{ij}'$, while with probability $1/(2e)$ the two 
processes jump exactly once by some $\pi$ and $\pi^*$ picked from the
joint distribution described above.  

Define $l(u) = v$ if this last possibility occurs (jumps of $\pi$ and 
$\pi^*$) and if furthermore, $\pi = \pi_{i(a,b),j(a,b)}$ and 
$\pi^* = \pi_{i^*(a,b),j^*(a,b)}$ for some $a,b \leq n/(6\sqrt{k})$.
This is of course measurable with respect to events until time $T(u,v)+1$,
and when it occurs, $W(T(u,v)+1) = W(T(u,v)) \cup \{ x_s \})$,
with $x_s = x_{s(u,v)}$ being the witnessing card for the stopping
time $T(u,v)$.  In this case, $T_u$ is defined to equal $T (u,v) + 1$
and the inductive statement~(\ref{W size}) is verified.
On the other hand, if $l(u) > v$, then $W(T(u,v)+1)
= W(T(u,v))$, since either the shuffles evolved in parallel or 
else~(\ref{stayeq}) guarantees that no card $x_{s'}$ other than $x_s$ 
was moved by either shuffle.  Thus again,~(\ref{W size})~is verified.
In either case (parallel shuffles or no card $x_{s'}$ other than $x_s$
moved by either shuffle), it is clear that the statement $\sigma_t (x_r)
\neq \tau_t (x_{r'})$ is preserved for all $r \neq r'$.

A consequence of~(\ref{W size}) is that all $k$ cards are
coupled by time $T_k$.  Thus to prove the theorem it suffices
to find a constant $c$ for which
\begin{equation} \label{tobe}
\P [T_k > cnk^3 (\ln(k))^2] < 1/2 .
\end{equation}
Let $\F(t)$ denote the $\sigma$-field of events up to time $t$.
We begin by showing that $\E T_0 < cn \ln (k)$.
Using Lemma~\ref{k=1} for $t = c_0 n \ln (k)$, with 
$c_0 > 3 c / \ln 2$, gives 
$$||\S^r - U||_1 < {1 \over k^3} .$$ 
Then for this $t$, $\P (\sigma_t (x_s) = \tau_t (x_{s'})) 
\leq 1/k^3 + 1/n^2$ for each fixed $s ,s' \leq k$ 
and summing gives a probability of at most $1/k + k^2/n^2$ that 
some $\sigma_t (x_s) = \tau_t (x_{s'})$.  Since $4 \leq k < \sqrt{n}$
in any nontrivial case, this probability is bounded above by $1/2$.
Repeating this argument at times that are multiples of $t$ shows
$T_0$ to be stochastically dominated by $t$ times a geometric
of mean two, proving that $\E T_0 < c n \ln (k)$ for an appropriate $c$.

Next, we establish that 
$$\begin{array}{rcrcl}
(i) & ~~~ & \E (T(u,v+1) - (T(u,v)+1) \| \F (T(u,v)+1)) & \leq 
   &  cn {k^2 \ln (k) \over k+1-u} \\
(ii) & ~~~ & \E (T(u,1) - T_{u-1}\| \F (T_{u-1})) & \leq 
   &  cn { k^2  \ln (k) \over k+1-u} \end{array} . $$
By Lemma~\ref{k=3}, choose $r = cn \ln (k)$ so that 
$||\S^r - U||_3 < 1/(400k^5)$.  Write $B$
for the region $[1,n/(3\sqrt{k})] \times [1,n/(3\sqrt{k})]$ and
write $C$ for the region $[1,n/(2\sqrt{k})] \times [1,n/(2\sqrt{k})]$.
Pick any $s \notin W(u)$ and let $y = \sigma_{T(u,v)+1} (x_s)$ and
$z = \tau_{T(u,v)+1} (x_s)$.  The set $Q$ of permutations $\pi$
for which $\pi (y) \in B$ and $\pi (z) \in B$ has probability
$$U (Q) = {1 \over 9k} \left ( {1 \over 9k} - {1 \over n^2} \right )
   \geq {1 \over 100 k^2} $$
under the uniform distribution.  The permutations $\sigma_{T(u,v)+1+r}
(\sigma_{T(u,v)+1}^{-1})$ and $\tau_{T(u,v)+1+r}
(\tau_{T(u,v)+1}^{-1})$ are equal and their conditional distribution
given $\F (T(u,v)+1)$ is the distribution of $\S^r$.  Since $r$ is
chosen to make $||\S^r - U||_2 \leq ||\S^r - U||_3 < 1/(400k^5)$,
it follows that
\begin{equation} \label{eq n1}
\P (\sigma_{T(u,v)+1+r} (x_s) \in B \mbox{ and } \tau_{T(u,v)+1+r} 
   (x_s) \in B \| \F (T(u,v)+1)) \geq {1 \over 100 k^2} - 
   {1 \over 400k^5} \; .
\end{equation} 
For $w \neq y,z$, 
$$U \{ \pi : \pi (y) \in B , \pi (z) \in B \mbox{ and } \pi (w) \in C \}
   \leq {1 \over 324 k^2 } .$$
Setting $w = \sigma_{T(u,v)+1} (x_{s'})$ for some $s' \neq s$ and
using $|| \S^r - U||_3 \leq 1/(400 k^5)$ again yields
\begin{eqnarray} \label{eq n2}
&& \P (\sigma_{T(u,v)+1+r} (x_s) \in B \mbox{ and } \tau_{T(u,v)+1+r} 
   (x_s) \in B \mbox{ and } \sigma_{T(u,v)+1+r} (x_s') \in C \|
    \F (T(u,v)+1))  \\
& \leq & {1 \over 324k^3}+{1 \over 400 k^5} .  \nonumber
\end{eqnarray}
If we instead let $w = \tau_{T(u,v)+1} (x_{s'})$, we see that the same is 
true with $\sigma_{T(u,v)+1+r} (x_s') \in C$ replaced by 
$\tau_{T(u,v)+1+r} (x_s') \in C$.  Let $G(u,v,s)$ be the event
that $\sigma_{T(u,v)+1+r} (x_s) \in B$, that $\tau_{T(u,v)+1+r} 
(x_s) \in B$, and that for all $s' \neq s$, $\sigma_{T(u,v)+1+r} (x_s'), 
\tau_{T(u,v)+1+r} (x_s') \notin C$.  Then summing~(\ref{eq n2}) over 
$s' \neq s$, doubling, and subtracting from~(\ref{eq n1}), gives
$$ \P (G(u,v,s) \| \F(T(u,v)+1)) \geq {1 \over 100 k^2} - 
   {1 \over 400 k^5} -  2k \left ( {1 \over 324 k^3} + {1 \over 400 k^5}
    \right ) \geq {1 \over 400 k^2} \, ,$$
since $k \geq 4$.  
The events $G(u,v,s)$ are disjoint as $s$ varies.  Recalling that $T(u,v+1)$
has been reached when $G(u,v,s)$ occurs for some $s \notin W(T_{u-1})$
and summing over such $s$ gives
$$\P (T(u,v+1) \leq T(u,v) + 1 + r \| \F(T(u,v)+1)) \geq 
    {k+1-u \over 400 k^2} .$$
Comparing to another geometric random variable, recalling the
value of $r$ and rolling all constants into one gives
$$\E (T(u,v+1)- T(u,v) - 1 \| \F (T(u,v)+1)) \leq cn {k^2 \ln (k)
   \over k+1-u} .$$
This establishes $(i)$ above, the argument for $(ii)$ being identical.

By construction, $\P (T_u = T(u,v) + 1 \| \F(T(u,v))) \geq {1 \over 36k}$
on the event $l(u) \geq v$.  This implies $\E l(u) \leq 36k$.  Thus,
setting $T(u,0) = T_{u-1}$, 
\begin{eqnarray*}
\E (T_u - T_{u-1}) & = & \E \sum_{v=1}^{l(u)} (T(u,v) -  T(u,v-1)) 
   \\[2ex]
& = & \E \sum_{v=1}^\infty \one_{l(u) \geq v} (T(u,v) -  T(u,v-1)) \\[2ex]
& = & \E \left [ \sum_{v=1}^\infty \one_{l(u) \geq v} \E (T(u,v) -  
    T(u,v-1)-1 \| \F(T(u,v)+1)) \right ] + \E l(u) \\[2ex]
& \leq & \E \left [ \sum_{v=1}^\infty \one_{l(u) \geq v} cn {k^2 \ln (k)
   \over k+1-u}  \right ] \\[2ex]
& = & \E l(u) cn {k^2 \ln (k) \over k+1-u}.
\end{eqnarray*}
Summing over $u$ gives $\E (T_k - T_0) \leq cn k^3 (\ln (k))^2$,
and using the earlier estimate on $\E T_0$ shows that
$\E T_k \leq cn k^3 (\ln(k))^2$.  Since $T_k$ is positive, this
implies~(\ref{tobe}), which proves Theorem~\ref{k-sets}.    $\Cox$

\section{Proof of Theorem \protect{\ref{tv time}}}

The proof of the nontrivial part of Theorem~\ref{tv time}, 
namely the upper bound, is gotten by analyzing the eigenvalues
of the random walk on $S_{n^2}$ whose steps have distribution $\S$.
To abbreviate the terminology, say the eigenvalues of a probability 
distribution $\P$ are the eigenvalues of its random walk, and
if $\P$ is uniform on some set $A$, call these also the eigenvalues
of $A$.  

The eigenvalue analysis is done in three steps.  Define another
shuffle $\R$ which chooses a three-cycle uniformly from among
all $2 {n^2 \choose 3}$ three-cycles at total rate one.  (A three-cycle 
permutes three cards cyclically and leaves the remaining $n^2-3$ cards
untouched.)  The first step, Lemma~\ref{compare} below,
compares the eigenvalues of $\S$ with the eigenvalues
of $\R$,  This  relies on a lemma from \cite{DS}, Lemma~\ref{diaconis} 
below, which bounds
the eigenvalues of one shuffle in terms of the eigenvalues of a
second, more tractable, shuffle when
the permutations in the second shuffle are explicitly written as 
products of permutations in the first shuffle.  The second
step is to compute the eigenvalues of $\R$.  This is done via 
the representation theory of the symmetric group, and can be
read off from known results in \cite{In}.  Finally, the information
about the eigenvalues of $\S$ is used to get an upper bound on the 
difference between $\S^t$ and $U$ in total variation, and hence
on the time to randomization.  This argument closely parallels
the proof of Theorem~5 in \cite[ch. 3]{Di}, which does an analogous
computation but for transpositions instead of three-cycles.

\begin{lem}[Diaconis 1992] \label{diaconis}
Let $A_1 , A_2 \subseteq S_n$ be sets of permutations that
generate $S_n$ and are symmetric, i.e. $\pi \in A_i$
if and only if $\pi^{-1} \in A_i$.  For each $\pi \in A_2$, pick
a way of writing $\pi$ as a product of elements of $A_1$;
let $N(\sigma , \pi)$ denote the number of times $\sigma$
appears in this product and let $|\pi|$ denote the number of
factors in the product.  This defines a constant
$$B =  { |A_1| \over |A_2| } \max_{\sigma \in A_1} \sum_{\pi
  \in A_2} |\pi| N(\sigma , \pi ) .$$
Let $\S_i$ be the uniform distribution on $A_i$.  
Choose any subspace $V \subseteq {\bf C}^{S_n}$ which is 
invariant for the right regular representation of $S_n$ and 
let $\lambda_1 \geq \lambda_2 \geq \cdots \geq \lambda_k$ be the
eigenvalues of $\S_2$ on the subspace $V$ in descending order,
counted with proper multiplicity.  Writing the eigenvalues of $\S_1$ 
on the subspace $V$ as $\lambda_1' \geq \lambda_2' 
\geq \cdots \geq \lambda_k'$, the relation 
\begin{equation} \label{comp bound}
1 - \lambda_i  \leq B (1 - \lambda_i')
\end{equation}
holds for $i = 1 , \ldots , k$.
$\Cox$
\end{lem}

\noindent{Proof:}  Let $\cE$ be the Dirichlet form for $\S_2$, namely
the symmetric, positive definite form on ${\bf C}^{S_n}$ defined
by $\cE (f,f) = \; < (I - \S_2) (f) , f >$, where $\S_2 (f) (z) = 
|A_2|^{-1} \sum_{x \in A_2} f(zx)$ and $<,>$ is the usual inner product.
Let $\cE'$ be the Dirichlet form for $\S_1$.  Then Theorem~1 
of~\cite{DS} shows that 
$$\cE \leq B \cE' .$$
Lemma~4 of~\cite{DS}
then implies~(\ref{comp bound}) when $V$ is all of ${\bf C}^{S_n}$.
If $V$ is not the whole space, then observe that $V$ has an orthogonal
complement $V^\perp$ which is also an invariant subspace.  Thus 
the Dirichlet forms $\cE$ and $\cE'$ decompose into the direct sums
of forms on $V$ and $V^\perp$.  The relation $\cE \leq B \cE' $ 
must then hold on $V$, and the proof is again finished by
Lemma~4 of~\cite{DS}.    $\Cox$

\begin{lem} \label{compare}
Let $\lambda_1 \geq \lambda_2 \geq \cdots \geq \lambda_{n!-2}$ 
be all the eigenvalues of the shuffle $\R$ except for the two 
eigenvalues of $+1$ which occur on the one-dimensional invariant 
subspaces $V_+ = \{ f : f(x) = f(y) \mbox{ for all } x,y \}$
and $V_- = \{ f : f(x) sign(x) = f(y) sign (y) \mbox{ for all } 
x , y \}$.  Let $\lambda_1' \geq \lambda_2', \geq 
\cdots \geq \lambda_{n!-2}'$ be the eigenvalues of $\S_0$ on the
space $V_\perp = \{ f : \sum f(x) = \sum f(x) sign (x) = 0 \}$
which is the orthogonal complement of $(V_+ \oplus V_-)$.
There is a constant $c$ such that for all $i \leq n! - 2$, 
$$ (1 - \lambda_i) \leq c (1 - \lambda_i').$$
The same holds when $\S_0$ is replaced by $\S$.
\end{lem}

\noindent{Proof:}  We first handle the case of $\S_0$.
To apply Lemma~\ref{diaconis}, let $A_1$ be all
the $\pi_{ij}$ and let $A_2$ be all the three-cycles.  Picking ways
to write elements of $A_2$ as products of elements of $A_1$ requires 
several steps.  Let $A_3 \subseteq A_2$ be the three-cycles that permute 
three array elements $(i_r , j_r) : r = 1,2,3$ for which the coordinates
$i_r$ are distinct from each other and the coordinates $j_r$ are
distinct from each other.  For $n \geq i,j \geq 3$, let $X_{ij}$ and $Y_{ij}$
be the following  product of elements of $A_1$ (commas are 
introduced for clarity and the notation for products is left-to-right,
so that $\pi \sigma$ means first do $\pi$ then $\sigma$):
\begin{eqnarray*}
X_{ij} & \deq & \pi_{i,j} \pi_{i-1,j} \pi_{i-2,j} \pi_{i-1,j} \\
Y_{ij} & \deq & X_{i,j} X_{i,j-1} X_{i,j-2} X_{i,j-1} .
\end{eqnarray*}
For $n \geq i \geq 3 > j$, let $X_{ij}$ be defined as above and let
$Y_{ij} = X_{ij}$.  For $n \geq j \geq 3 > i$, let $X_{ij} = \pi_{ij}$
and let $Y_{ij} =  X_{i,j} X_{i,j-1} X_{i,j-2} X_{i,j-1}$ as before.
Finally, if $3 \geq i,j$, let $X_{ij} = Y_{ij} = \pi_{ij}$.  

Claim: $Y_{ij}$ is the permutation that transposes the $i,j$-element
of the array with the top element $T$, and in addition, if $i,j \geq 2$,
transposes the $i,1$-element with the $1,j$-element.  The proof
of this is omitted, being a case by case verification; the figure
illustrates the case $i=j=5$.  \\[2ex]
$$ X_{55} \; : \; \begin{array}{ccccc} 11 & 12 & 13 & 14 & 15 \\ 
   21 & 22 & 23 & 24 & 25 \\ 31 & 32 & 33 & 34 & 35 \\ 
   41 & 42 & 43 & 44 & 45 \\ 51 & 52 & 53 & 54 & 55 \end{array}
   \longrightarrow_{(\pi_{55})}  
   \begin{array}{ccccc} 55 & 54  & 53 & 52 & 51 \\ 
   45 & 44 & 43 & 42 & 41 \\ 35 & 34 & 33 & 32 & 31 \\ 
   25 & 24 & 23 & 22 & 21 \\ 15 & 14 & 13 & 12 & 11 \end{array}
   \longrightarrow_{(\pi_{45})}  
   \begin{array}{ccccc}   21 & 22 & 23 & 24 & 25 \\ 
   31 & 32 & 33 & 34 & 35 \\ 41 & 42 & 43 & 44 & 45 \\ 
   51 & 52 & 53 & 54 & 55 \\ 15 & 14 & 13 & 12 & 11  \end{array} 
   \hspace{.5in}$$  
$$ \hspace{1.5in} \longrightarrow_{(\pi_{35})}  
   \begin{array}{ccccc} 45 & 44 & 43 & 42 & 41 \\ 
   35 & 34 & 33 & 32 & 31 \\ 25 & 24 & 23 & 22 & 21 \\ 
   51 & 52 & 53 & 54 & 55 \\ 15 & 14 & 13 & 12 & 11 \end{array}
   \longrightarrow_{(\pi_{45})}  
   \begin{array}{ccccc} 55 & 54  & 53 & 52 & 51 \\ 
   21 & 22 & 23 & 24 & 25 \\ 31 & 32 & 33 & 34 & 35 \\ 
   41 & 42 & 43 & 44 & 45 \\ 15 & 14 & 13 & 12 & 11 \end{array}$$
\vspace{.4in}
$$ Y_{55} \; : \;  \begin{array}{ccccc} 11 & 12 & 13 & 14 & 15 \\ 
   21 & 22 & 23 & 24 & 25 \\ 31 & 32 & 33 & 34 & 35 \\ 
   41 & 42 & 43 & 44 & 45 \\ 51 & 52 & 53 & 54 & 55 \end{array}
   \longrightarrow_{(X_{55})} 
   \begin{array}{ccccc} 55 & 54  & 53 & 52 & 51 \\ 
   21 & 22 & 23 & 24 & 25 \\ 31 & 32 & 33 & 34 & 35 \\ 
   41 & 42 & 43 & 44 & 45 \\ 15 & 14 & 13 & 12 & 11 \end{array}
   \longrightarrow_{(X_{54})} 
   \begin{array}{ccccc} 12 & 13 & 14 & 15 & 51 \\ 
   21 & 22 & 23 & 24 & 25 \\ 31 & 32 & 33 & 34 & 35 \\ 
   41 & 42 & 43 & 44 & 45 \\ 51 & 53 & 54 & 55 & 11 \end{array} 
   \hspace{.5in}$$
$$ \hspace{1.5in} \longrightarrow_{(X_{53})} 
   \begin{array}{ccccc} 54  & 53 & 52 & 15 & 51 \\ 
   21 & 22 & 23 & 24 & 25 \\ 31 & 32 & 33 & 34 & 35 \\ 
   41 & 42 & 43 & 44 & 45 \\ 14 & 13 & 12 & 55 & 11 \end{array} 
   \longrightarrow_{(X_{54})} 
   \begin{array}{ccccc} 55 & 12 & 13 & 14 & 51 \\ 
   21 & 22 & 23 & 24 & 25 \\ 31 & 32 & 33 & 34 & 35 \\ 
   41 & 42 & 43 & 44 & 45 \\ 15 & 52 & 53 & 54 & 11 \end{array}$$

Next, for pairs $(i_1 , j_1) , (i_2 , j_2)$ both unequal to $T$
and satisfying $i_1 \neq i_2$ and $j_1 \neq j_2$, let 
$$Z_{i_1 , j_1 , i_2 , j_2} = Y_{i_i  , j_1} Y_{i_2  , j_2} Y_{i_1 , j_1} 
   Y_{i_2 , j_2} . $$
It is easy to see that $Z_{i_1 , j_1 , i_2 , j_2}$ is the three-cycle 
permuting $T$, the $i_2 , j_2$-element and the $i_1 , j_1$-element.
Finally, for $i_1 , j_1 , i_2 , j_2 , i_3 , j_3$ with none of
the $i_r$'s equal to another, none of the $j_r$'s equal to another and
no pair $(i_r , j_r)$ equal to $(1,1)$, let 
$$W_{i_1 , j_1 , i_2 , j_2 , i_3 , j_3} = Z_{i_1 , j_1 , i_2 , j_2}
Z_{i_2 , j_2 , i_3 , j_3} .$$
Then $W_{i_1 , j_1 , i_2 , j_2 , i_3 , j_3}$ cyclically permutes
the $i_3 , j_3$-element, the $i_2 , j_2$-element and the 
$i_1 , j_1$-element.  If $\pi \in A_3$ is a three-cycle that
permutes three array elements $(i_3 , j_3) , (i_2 , j_2)$ and
$(i_1 , j_1)$ with $i_r , j_r \geq 2$, pick the decomposition of
$\pi$ into elements of $A_1$ according to the construction of
$W_{i_1 , j_1 , i_2 , j_2 , i_3 , j_3}$; if one of the pairs $(i_r ,
j_r)$ is equal to $(1,1)$, then use the appropriate $Z$ instead
of $W$.  In the obvious notation, 
$|\pi| = | Z_{i_1 , j_1 , i_2 , j_2} | +  | Z_{i_2 , j_2 , i_3 , j_3} |
\leq 128$.  Furthermore, for any $\sigma = \pi_{ij} \in A_1$, the number
of $\pi \in A_3$ for which $N(\sigma , \pi) > 0$ is at most
$27 n^4$, since one of the pairs $(i_r , j_r)$ must satisfy 
$i \leq i_r \leq i+2$ and $j \leq j_r \leq j+2$.  Thus
$$\sum_{\pi \in A_3} N(\sigma , \pi) \leq 32 \cdot 27 n^4$$
for any $\sigma \in A_1$. 

For $\pi \in A_2 \setminus A_3$, decompose it into a product
of elements of $A_1$ as follows.  If $\pi$ permutes the $i_r,j_r$-elements
for $r = 1,2,3$, choose $(u_1, v_1)$ and $( u_2 , v_2)$ from among
the set $\{ (x,y) :  \exists r \mbox{ with } |x - i_r| + |y - j_r| \leq 6 \}$
in such a way that each $u_s$ is distinct from each $i_r$, each
$v_s$ is distinct from each $j_r$, and $u_1 \neq u_2$ and
$v_1 \neq v_2$.  Writing $a,b,c,d,e$ for
$(i_1 , j_1) , (i_2 , j_2) , (i_3 , j_3) , (u_1 , v_1)$, and $(u_2, v_2)$
respectively, decompose $\pi$ as
$$ \pi = W_{ade} W_{bde} W_{cde} W_{ade} W_{bde} .$$
It is easy to check that this does indeed give $\pi$ and that
for $\pi \in A_2 \setminus A_3$, the decomposition satisfies
$|\pi| \leq 640$.  Furthermore, the number of $\pi \in A_2 \setminus A_3$
for which $N(\pi_{ij} , \pi) > 0$ is bounded by the number
of ways of choosing three array elements in such a way
that some two are in the same row or column and one is
within a distance 6 of $(i,j)$ in the taxicab metric.  This
is at most $c n^4$ for some constant $c$.

Applying Lemma~\ref{diaconis} with $n > 10$ now gives
$1 - \lambda_i  \leq B (1 - \lambda_i')$ where 
\begin{eqnarray*}
B  & = &  { |A_1| \over |A_2| } \max_{\sigma \in A_1} \sum_{\pi 
  \in A_2} |\pi| N(\sigma , \pi ) \\[2ex]
& \leq & {n^2 \over 2 {n^2 \choose 3}} \left ( \max_{\sigma \in A_1} \sum_{\pi 
  \in A_3} |\pi| N(\sigma , \pi ) + \max_{\sigma \in A_1} \sum_{\pi 
  \in A_2 \setminus A_3} |\pi| N(\sigma , \pi ) \right )  \\[2ex]
& \leq & 3.1 n^{-4} (128 \cdot (128 \cdot 27 n^4) + 640 \cdot (640 c n^4))  
\end{eqnarray*}
which is bounded by some constant, proving the lemma for
$\S_0$.  For $\S$, use the same decompositions, losing
a factor of two in $|A_1| / |A_2|$.    $\Cox$

It has been shown that the eigenvalues of $\S$ are bounded 
in terms of the eigenvalues of $\R$; the computation of these
latter uses a combinatorial formula from~\cite{In}.
Let $\rho$ be any irreducible matrix representation of $S_n$.
Since the measure $\R$ is uniform on conjugacy classes,
the matrix $\R (\rho) \deq \E_{\R} (\rho)$ will be a
constant multiple of the identity, the constant being
$\chi_\rho (\tau) / d(\rho)$, where $\chi_\rho$ is the
character of the representation $\rho$ and $\tau$ 
is any element of the
conjugacy class, in other words, any three-cycle.  This
gives $d(\rho)$ eigenvalues equal to $\chi_\rho (\tau) / d(\rho)$
in the irreducible representation $\rho$, and since
this representation appears with multiplicity $d(\rho)$,
the shuffle $\R$ will have this eigenvalue with 
multiplicity $d(\rho)^2$.  Ingram's formula for the
characters of the irreducible representations of $S_n$ evaluated 
at a three-cycle yields the following upper bounds:

\begin{lem} \label{char}
Let $\rho$ be the irreducible representation of $S_n$ 
corresponding to the partition $t = (t_1 \geq t_2 \geq \cdots)$
of $n$.  Then the character of $\rho$ evaluated at a three-cycle
is given by 
\begin{equation} \label{Ingram}
r (\rho) \deq \chi_\rho (\tau) / d(\rho) = { 3 \sum_{i,j} (i-j)^2 
  \over n(n-1)(n-2)} \, - \, {3 \over 2(n-2)} , 
\end{equation}
where the sum is over all $(i,j)$ such that $t_i \geq j$,
or in other words over all squares of the Young tableau
for the partition $t$.  It follows from this that
$$ r (\rho) \leq 1 - { 3 (t_1 - 1) (n-t_1) \over (n-1)(n-2) }
  \mbox{   when } t_1 \geq n/2 $$
and 
$$ r (\rho) \leq \max \{ t_1 - 1, t_1'- 1 \} / (n-2)
  \mbox{   when } t_1 , t_1' \leq n/2 , $$
where $t_1' = \max \{ i : t_i > 0 \}$ is the first element of
the partition dual to $t$.  
\end{lem}

\noindent{Proof:}  The formula~(\ref{Ingram}) is taken directly
from \cite[(5.2)]{In}, where the term $a(a+1)(2a+1)$ is replaced
by $6 \sum_{i=1}^a i^2$ and the typographical error (a misplaced
parenthesis) is corrected.  For fixed $t_1 \geq n/2$, the sum is 
maximized by letting $t_2 = \cdots = t_{n+1-t_1} = 1$ and $t_i = 0$ for
$i > n+1-t_1$.  For the trivial representation, $t = n,0,0,\ldots$
and $r = 1$.  Comparing~(\ref{Ingram}) for the trivial representation
and a nontrivial representation $\rho$ gives 
\begin{eqnarray*} 
1 - r(\rho) & \geq & {3 \over n(n-1)(n-2)} \left [
    \sum_{k=1}^{n - t_1} ((t_1 - 1 + k)^2 - k^2) \right ] \\[2ex]
& = & {3 \over n(n-1)(n-2)} \left [
    \sum_{k=1}^{n-t_1} (t_1 - 1)^2 + 2k (t_1 - 1) \right ] \\[2ex]
& = & {3 \over n(n-1)(n-2)} \left [ 
    (n - t_1) (t_1 - 1)^2 - (n - t_1) (n - t_1 + 1) (t_1 - 1) \right ] \\[2ex]
& = &  {3 \over (n-1)(n-2)} \left [ (n- t_1) (t_1 - 1) \right ] .
\end{eqnarray*}
On the other hand, when $t_1 , t_1' \leq n/2$, then let $t_0
= \max \{ t_1 , t_1' \}$.  Ignore the subtracted term in~(\ref{Ingram})  
to get 
$$ r(\rho) < { 3 \sum_{i,j} (i-j)^2 \over n(n-1)(n-2)} .$$
Partition the $n$ pairs $(i,j)$ according to the value of $i$
and observe that for any $i$, the average of the summands
with that particular value of $i$ is
\begin{eqnarray*}
t_i^{-1} \sum_{j = 1}^{t_i} (j - i)^2 & \leq & t_0^{-1} \sum_{j=1}^{t_0}
   (j-1)^2  \\[2ex]
& = & (t_0 - 1) (2t_0 - 1) / 6 .
\end{eqnarray*}
This is then an upper bound for the average of all the summands; the
sum is precisely $n$ times the average, yielding
$$ r(\rho) < {(t_0 - 1) (2 t_0 - 1) \over 2 (n-1) (n-2) } 
   \leq {t_0 - 1 \over 2 (n-2)} \; .$$ 
$\Cox$

The bound~(\ref{UBL}) below on the time to randomization 
for the shuffle $\S$ in terms of its eigenvalues is based on
the Upper Bound Lemma (3b.1) from~\cite{Di}; the evaluation
of~(\ref{UBL}) is based on the analogous computation for
random transpositions on pages 41 - 42 of~\cite{Di}.  Accordingly,
some details are omitted here.
 
\noindent{Proof} of Theorem~\ref{tv time}:  Let the eigenvalues of 
$\R$ and $\S$ be denoted respectively by $\lambda_i$ and
$\lambda_i'$, listed in the following order: $\lambda_1 = 
\lambda_1' = 1$ are the eigenvalues on $V_+$; $\lambda_2, \lambda_2'$ 
are the eigenvalues on $V_-$, with $\lambda_2 = 1 > \lambda_2'$; 
$\lambda_3 \geq \cdots \geq \lambda_{n!}$ and $\lambda_3' \geq 
\cdots \geq \lambda_{n!}'$ are the eigenvalues on $V_\perp$.
Using the constant $c$ from Lemma~\ref{compare}
and Lemma~3B.1 of~\cite{Di} gives
\begin{eqnarray}
4 |\S^{ct} - U|^2 & \leq n! & \sum |\S^{ct} (\pi ) - U(\pi)|^2 
    \nonumber \\[2ex]
& = & \sum_{i \geq 2} e^{-2ct(1-\lambda_i')} \nonumber \\[2ex]
& \leq &  e^{-2ct (1-\lambda_2')} + \sum_{i \geq 3} e^{-2t(1-\lambda_i)} 
    \nonumber \\[2ex]
& = & e^{-2ct (1-\lambda_2')} + {\sum_{\rho}}^* d(\rho)^2 \exp 
    [-2t (1-r(\rho))],  \label{UBL}
\end{eqnarray} 
where $\sum^*$ denotes a sum is over representations $\rho$ other 
than the trivial representation and the alternating representation.

We now bound~(\ref{UBL}) using Lemma~\ref{char}.  First dispose of
the $e^{-2ct (1 - \lambda_2')}$ term.  Since the alternating
character is $\sum \sgn (\sigma) \S (\sigma)$ and the sign
of $\pi_{ij}$ is negative when (among other cases) $i$ is odd and
$j \equiv 2$ mod 4, the alternating character is at most $3/4$, and
$$e^{-2ct (1-\lambda_2')} \leq e^{-ct/2} .$$

For the remaining sum, observe that if $\rho$ and $\rho'$ correspond
to dual partitions $t,t'$ then $d(\rho) = d(\rho')$ and $r(\rho) =
r(\rho')$.  Since the trivial and alternating partitions are dual, this gives
$$  {\sum_{\rho}}^* d(\rho)^2 \exp [-2t (1-r(\rho))] \leq
    2{\sum_{\rho}}^{**} d(\rho)^2 \exp [-2t (1-r(\rho))]$$
where $\sum^{**}$ is over nontrivial partitions with $t_1 \geq t_1'$.  
Note that for $t \geq n/2$,
\begin{eqnarray*}
1 - { 3 (t - 1) (n - t) \over n (n-1) } & = & 1 - {3 (t-1) \over n-2}
   \left ( 1 - {t - 1 \over n-1} \right ) \\[2ex]
& \leq & 1 - {3 \over 2} \, (1 - {t-1 \over n-1}) \\[2ex]
& \leq & {t-1 \over n-2}
\end{eqnarray*}
and thus for any $\alpha \in (0,1/2)$, the above expression 
involving $\sum^{**}$ is at most
$$ 2 \sum_{\rho : t_1 \geq (1-\alpha)n}^{**} d(\rho)^2 \exp [-2t
    { 3 (t_1 - 1) (n-t_1) \over (n-1)(n-2) }]
   + 2 \sum_{\rho : t_1 < (1-\alpha)n} d(\rho)^2 \exp [-2t
    {t_1 - 1 \over n-2} ] .$$
Diaconis now shows~\cite[proof of Theorem 5, page 42]{Di} that
$\alpha \in (0,1/4)$ may be chosen so that when $t > (1/2) n \ln (n)
+ kn$, both sums together are less than $a e^{-2k}$ for some universal 
constant $a$.  This shows that $|\S^{ct} - U|$ goes to zero when $t =
(.5 + \ee) n \ln (n)$, proving Theorem~\ref{tv time}.   $\Cox$

\renewcommand{\baselinestretch}{1.0}\large

\noindent{Department of Mathematics} \\
University of Wisconsin-Madison \\
Van Vleck Hall \\
480 Lincoln Drive \\
Madison, WI 53706

\end{document}